\documentclass[10pt]{article}
\usepackage[cp850]{inputenc}
\usepackage[T1]{fontenc}\usepackage{amsfonts}
\usepackage{graphicx}
\usepackage{amsmath}
\def\N{\mathbb{N}}\def\E{\mathbb{E}}
\def\d{\frac{1}{2}\, }
\author{G\'erard  Letac\thanks{Institut de Math\'ematiques de Toulouse, Universit\'e Paul Sabatier, 118 route de Narbonne 31062 Toulouse, \texttt{gerard.letac@math.univ-toulouse.fr}}}
\date{\today}

\title{Is the Sibuya distribution  a progeny? }
\begin{document}  \maketitle

\begin{abstract} For $0<a<1$ the Sibuya distribution $s_a$ is concentrated on the set $\mathbb{N}^+$ of positive integers and is defined by  the generating function $\sum_{n=1}^{\infty}s_a(n)z^n=1-(1-z)^a.$ A distribution $q$ on $\mathbb{N}^+$  is called a progeny if there exists  a  Galton-Watson process 
 $(Z_n)_{n\geq 0}$ 
such that $Z_0=1$, such that  $\E(Z_1)\leq 1$ and such that $q$ is the distribution of $\sum _{n=0}^{\infty}Z_n. $ The paper proves that  $s_a$ is a progeny if and only if $\d\leq a<1.$ The point is to find the values of $b=1/a$ such that the power series expansion of $u(1-(1-u)^b)^{-1}$ has non negative coefficients. The proof is not short, but elementary.

\end{abstract}

\section{Introduction} \subsection{Branching processes}Recall that if $p$ is a distribution on the set $\N$ of non negative integers with generating function $$f_p(z)=\sum_{n=0}^{\infty}p_nz^n$$ then  a branching process $(Z_n)_{n=0}^{\infty}$ (also called Galton-Watson process) governed by $p$ 
 is the Markov chain defined on  $\N$ by $Z_0=1$ and $$Z_{n+1}=\sum_{k=1}^{Z_n}X_{n,k}$$ where $(X_{n,k})_{n,k\geq 0}$ are iid random variables with distribution $p.$
It is easily seen that
\begin{equation}\label{HARRIS}\E(z^{Z_n})=f_p^{(n)}(z)=f_p\circ\cdots\circ f_p(z)\ \  n\ \mathrm{times}.\end{equation} 
Harris (1963) and Athreya and Ney (1972) are two classical treatises on the subject. 

\subsection{Sibuya distribution} In all the sequel we use the Pchhammer symbol $(x)_n=x(x+1)\ldots(x+n-1).$ The Sibuya distribution $ s_a$ with parameter $a\in (0,1)$ is a probability on positive integers $\mathbb{N}^+=\{1,2,\ldots\}$ such that if $S\sim s_a$ then \begin{equation}\label{SIB}f_{s_a}(z)=\E(z^S)=1-(1-z)^a.\end{equation}  It is easily seen that for $n\geq 1$
\begin{equation}\label{DENSITY}\Pr(S=n)=s_a(n)=\frac{1}{n!}a(1-a)(2-a)\ldots (n-1-a)=\frac{(-a)_n}{n!}.\end{equation} Although this law $s_a$ has certainly been considered before, traditionally one refers to Sibuya (1979) for its study. In the Sibuya paper, it appears as a particular case of what Sibuya calls a digamma law. It can be found in his formulas (16)  and (28) by taking $\gamma=-\alpha=a$ for getting our \eqref{DENSITY}. Surprizingly  \eqref{SIB} does not appear in his paper. 
If $S\sim s^{1/2}$  the distribution of $2S$ is well known as the law of first return to 0 of a simple random walk. An other probabilistic interpretation of $S\sim s_a$ is $S(w)=\min\{n: w\in A_n\}$ where $(A_n)_{n\geq 1}$ is a sequence of  independent events such that $\Pr(A_n)=a/n.$ To see this just compute $\Pr(S=n)$ by using $$\Pr(S>n)=(1-a)(1-\frac{a}{2})\ldots(1-\frac{a}{n})$$ and obtaining \eqref{DENSITY}. This idea is extended in Kozubowski and Podg\'orski (2018), where the authors rather consider $\Pr(A_n)=a/(n+k),$ when the real parameters $k$ and $a$ satisfy $0<a<k+1.$  For this generalization $s_{a,k}$ of $s_a=s_{a,0}$, the generating function of $S$ is expressed  with the hypergeometric function $_2F_1(k+1-a,1;k+1; z) $ with a simpler form (2.13 in the paper) when $k$ is a non negative integer: 
\begin{equation}\label{SIBk}f_{s_{a,k}}(z)=\E(z^S)= \frac{1}{z^kP_k(1)}
(P_k(z)-(1-z)^a),\ \ P_k(z)=\sum_{j=0}^k\frac{(-a)_j}{j!}z^j,\ \ P_k(1)=\frac{(1-a)_k}{k!},\end{equation} 
the value of $P_k(1)$ being computed by considering $\sum_{k=0}^{\infty}P_k(1)s^k.$ This article  of Kozubowski and Podg\'orsk contains numerous observations  and a rich set of references about $s_a.$

 However, the  most interesting feature  of the Sibuya law from the point of view of branching processes is the following semigroup property
$$f_{s_a}\circ f_{s_{a'}}=f_{s_{aa'}}$$ This implies that the branching process $(Z_n)_{n=0}^{\infty}$ is governed by $s_a$  then the law of $Z_n$ is quite explicit  and is $ s_{a^n}.$ This distribution is an important example of an imbeddable law  in a continuous semigroup for composition. See Grey (1975) for instance. There are some other variations of this distribution, sometimes  called informally also Sibuya distributions. One is  the mixing with an atom on zero  $(1-\lambda)\delta_0+\lambda\ s_a$ where $0<\lambda\leq 1,$ with generating function
 \begin{equation}\label{ATOME}f_{(1-\lambda)\delta_0+\lambda\ s_a}(u)=1-\lambda (1-u)^a\end{equation}
An other one is the natural exponential family extension of the Sibuya distribution, say $ s_a^{(\rho)},$ defined for $0<\rho\leq 1$ by its generating function
 \begin{equation}\label{EXP}f_{ s_a^{(\rho)}}(u)=\frac{1- (1-\rho u)^a}{1- (1-\rho)^a}.\end{equation}
However, in the present paper Sibuya  law will only mean $s_a$  for some $0<a<1.$ Last section will comment on \eqref{EXP}.

\subsection{Progeny} If $(Z_n)_{n=0}^{\infty}$ is a branching process governed by $p,$ a  classical fact is  that $\Pr(\exists n\ :\  Z_n=0)=1$ is and only if $m=\sum_{n=0}^{\infty}np_n\leq 1:$  consult Harris (1963) page 7. Under these circumstances the random variable 
$$S=\sum_{n=0}^{\infty}Z_n$$ is finite. Its distribution $q$ is called the progeny of $p$ and we have the following link between the generating functions of $p$ and $q$: for all $z$ such that $|z|\leq 1$ the following holds
\begin{equation}\label{PROGENY}f_q(z)=zf_p(f_q(z))
\end{equation} Progenies have been considered by many authors, and one can consult Harris (1963) page 32 for references. 
Since $Z_0=1$ the sum $S$ is concentrated on $\N^+.$ Given $p,$ the calculation of $q$, or  of $f_q,$ is not easy in general. It can be done by the Lagrange-B\"urmann formula (see Whittaker and Watson (1986) page 129 for instance). For calculating $q$, only the case where $f_p$ is a Moebius function is simple: a reference is Toulouse (1999) page 266. Surprizingly enough, given $f_q$ the calculation of $p$ is more feasible. Actually, the  function $z=g(u)=f_q^{(-1)}(u)$ valued in $[0,1]$ is well defined for $u\in [0,1]$ by $u=f_q(z)$ and \eqref{PROGENY} leads to 
\begin{equation}\label{PROGENY2}f_p(u)=\frac{u}{g(u)}.\end{equation}
Although the correspondence $p\mapsto q$ is one to one, clearly not all distributions $q$ on $\N^+$ can be the progeny of some $p.$ For instance $q$ cannot have a bounded support, except in the trivial case where $p_0=1.$ Another necessary condtion for $q$ to be a progeny is  $q_1=f'_q(0)>0$ . For if not, the reciprocal function $g$ of $f_q$ would not be analytic on $0.$  In general, given a probabilty $q$ on $\N^+$, no necessary and sufficient condition such that $q$ is a progeny is known. 
\subsection{When is the Sibuya distribution a progeny ?} Since the Sibuya distribution is concentrated on $\N^+$, has  an unbounded support and satisfies $q_1=a>0$, the natural question is the following: given $a\in (0,1),$ does there exists $p$ such that $q=s_a?$
The following proposition gives the answer and is the aim of the present paper:

\vspace{4mm}\noindent\textbf{Proposition 1.} The Sibuya distribution $s_a$ is a progeny if and only if $\d\leq a<1.$

\vspace{4mm}\noindent Section 2 gives the proof of Proposition 1.
Note that the question is a natural one, since the explicit calculation of $g$ is possible. Replacing $s_a$ by the generalization $s_{a,k}$  as described  in \eqref{SIBk} is clearly complicated for the calculation of $g.$ Also, one should not wonder whether  \eqref{ATOME} is a progeny, since it has an atom on zero. The case \eqref{EXP} is more interesting:  Section 3 comments on the new progeny  when the generating law $p$ is replaced by an element of the natural exponential family $p^{(r)}$ defined by its generating function $f_p(r)/f_p(r)$. These considerations are then applied  to \eqref{EXP}.

\section{Proof of  Proposition 1.} For simplicity we denote $b=1/a> 1.$ The calculation of the function $g$ appearing in \eqref{PROGENY2} is easy since if $u\in [0,1]$ the only solution $z\in [0,1]$ of the equation $u=1-(1-z)^{1/b}$ is $g(u)=1-(1-u)^b$. We have therefore to prove that the function $$u\mapsto h_b(u)=\frac{u}{1-(1-u)^b}$$ has a power series expansion with non negative coefficients if and only if $1\leq b\leq 2.$

In the sequel, we denote  by $B$ the set of  $b$'s such that $s_{1/b}$ is a progeny.

\vspace{4mm}\noindent\textbf{First step: 2 in in $B$, 3 is not in $B.$} 
We have $h_2(u)=\d\frac{1}{1-\d u}=\sum_{n=0}^{\infty}\frac{1}{2^{n+1}}u^n$ and we get back the well known fact that $s_{\d}$ is the progeny of a geometric distribution starting at 0. For $b=3$ we have
$$\ h_3(u)=\frac{1}{3}\times \frac{1}{1-u+\frac{u^2}{3}}=\frac{1}{3}\sum_{n=0}^{\infty}r^n\frac{\sin(n+1)\theta}{\sin \theta}u^n,$$
where $re^{\pm i\theta}=\d(3\pm i\sqrt{3})$  are the complex roots of the polynomial $1-u+\frac{u^2}{3}.$ Actually $r=\sqrt{3}$ and $\theta=\pm \pi/6.$ Clearly $\sin(n+1)\theta/\sin \theta \geq 0$ for all $n$ is impossible and therefore $3$ is not in $B.$

\vspace{4mm}\noindent\textbf{Second step: $b$ is in $B$ if $1<b<2.$ }

Denote
$$H(u)=\frac{b-1}{2}+\sum_{n=1}^{\infty}(b-1)(2-b)(3-b)\ldots(n+1-b)\frac{u^n}{(n+2)!}$$  Since $1<b<2$ all the coefficients of $H$ are positive. With the Pochhammer symbol we can write 
$$(1-u)^b=\frac{1}{(1-u)^{-b}}=\sum_{n=0}^{\infty}\frac{(-b)_n}{n!}u^n$$
As a consequence 
\begin{equation}\label{LOUIS}h_b(u)=\frac{1}{b}\times\frac{1}{1-uH(u)}=\frac{1}{b}\sum_{n=0}^{\infty}u^nH(u)^n.\end{equation} which implies that $b\in B.$

\vspace{4mm}\noindent\textbf{Third step: $b$ is not in $B$ if  $b>3.$}

Consider the numbers $(P_n)_{n\geq 2}$ defined by 
\begin{equation}\label{LESP}\frac{bu}{1-(1-u)^b}=1+\frac{b-1}{2}u+\sum_{n=2}^{\infty}P_nu^n\end{equation} A simple calculation shows that 
$$P_2=\frac{b^2-1}{6},\ P_3=\frac{b^2-1}{4},\ P_4=\frac{(19-b^2)(b^2-1)}{30},\ P_5=\frac{(9-b^2)(b^2-1)}{4}.$$ Therefore $ P_5<0$ if $b>3.$

\vspace{4mm}\noindent\textbf{Fourth step: $b$ is not in $B$ if  $2<b<3.$}

 This point is more difficult. Let us introduce the 
numbers $(p_n)_{n\geq 2}$ defined by 
\begin{eqnarray}\nonumber\frac{1-(1-u)^b}{bu}&=&\int_0^1(1-ux)^{b-1}dx\\&=&\label{LESp}1-\frac{b-1}{2}u+\sum_{n=2}^{\infty}\frac{(1-b)_n}{(n+1)!}u^n=1-\frac{b-1}{2}u+\sum_{n=2}^{\infty}p_nu^n\end{eqnarray} 
where the Pochhammer symbol is used:
\begin{equation}\label{pPOSIT}p_n=\frac{1}{(n+1)!}(b-1)(b-2)(3-b)\ldots (n-b)\end{equation}
Observe that from \eqref{pPOSIT} we have $p_n>0$  if $2<b<3$ since $n\geq 2.$ 
For simplication, we now denote 
$$v=\frac{b-1}{2}u,\ u=\frac{2}{b-1}v,\ A_n=P_n\frac{2^n}{(b-1)^n},\  a_n=p_n\frac{2^n}{(b-1)^n}.$$  With these notations, equalities \eqref{LESP} and \eqref{LESp}
become
\begin{eqnarray}\label{LESA}
\frac{2}{b-1}\times \frac{v}{1-\left(1-\frac{2}{b-1}v\right)^b}&=&1+v+\sum_{n=2}^{\infty}A_nv^n\\\label{LESa}
\frac{b-1}{2}\times \frac{1-\left(1-\frac{2}{b-1}v\right)^b}{v}&=&1-v+\sum_{n=2}^{\infty}a_nv^n
\end{eqnarray}
Since $2<b<3$ the radius of convergence of the power series on the right hand side of \eqref{LESa} is $(b-1)/2<1.$ Because $a_n>0$ this remark implies that 
\begin{equation}\label{DIV}\sum_{n=2}^{\infty}a_n=\infty.\end{equation}
Now we multiply the two right hand sides of \eqref{LESa} and \eqref{LESA}. This product is 1.  We obtain 
$$\sum_{n=2}^{\infty}(a_n+A_n)v^n+\sum_{n=2}^{\infty}(a_n-A_n)v^{n+1}+\left(\sum_{n=2}^{\infty}a_nv^n\right)\left(\sum_{n=2}^{\infty}A_nv^n\right)=v^2$$
$$\sum_{n=2}^{\infty}(a_n+A_n)v^n+\sum_{n=3}^{\infty}(a_{n-1}-A_{n-1})v^{n}+\sum_{n=4}^{\infty}\left(\sum_{k=2}^{n-2}A_{n-k}a_k\right)v^n=v^2$$
From this last equality, watching the coefficient on $v^n$ for $n\geq 4$ we get \begin{equation}\label{AHAH}a_n+a_{n-1}+\sum_{k=2}^{n-2}A_{n-k}a_k=A_{n-1}-A_n\end{equation} 
Now assume that $A_n\geq 0$ for all $n\geq 2,$  which is also assuming that $b\in B.$ Then \eqref{AHAH} implies 
$a_n\leq A_{n-1}-A_n.$ Summing up from $n=4$ to $N$ we get for all $N$
$$\sum_{n=4}^Na_n\leq A_3-A_N\leq A_3,$$ which contradicts \eqref{DIV}. Therefore  there exists at least one $n$ such that $A_n<0.$ The proposition is  proved. $\square$

\section{Natural exponential family and progeny}

\vspace{4mm}\noindent\textbf{Proposition 2.} Let $f_p$ be governing a branching process with mean $m\leq 1$ and with generating function $f_q$ for its progeny. Denote by $R\geq 1$ the radius of convergence of the power series $f_p.$ 
Consider $r\in (0,R)$ and ${f_{p^{(r)}}(z)=f_p(rz)/f_p(r)}.$  Suppose that that $m_r=rf_p'(r)/f_p(r)\leq 1.$ Then the progeny $q^{(\rho)}$ associated to $p^{(r)}$ is given by the following generating function
$$f_{q^{(\rho)}}(z)=f_q(\rho z)/f_q(\rho)\ \  \ \mathrm{with}\ \ r=f_q(\rho].$$

\vspace{4mm}\noindent\textbf{Comment.}
In other terms, if the branching process is governed by a distribution belonging to the natural exponential family generated by the probability $p$, the corresponding progeny has a distribution belonging to the natural exponential family  generated by $q$,  but with a new parameter.  

\vspace{4mm}\noindent\textbf{Proof.} \begin{eqnarray*}sf_{p^{(r)}}(f_{q^{(\rho)}}(z))&=&\frac{z}{f_p(r)}f_p(rf_{q^{(\rho)}}(z))=\frac{z}{f_p(r)}f_p\left(r\frac{f_q(\rho z)}{f_q(\rho)}\right)\\&=&\frac{\rho z}{\rho f_p(f_q(\rho))}f_p(f_q(\rho z))=\frac{f_q(\rho z)}{f_q(\rho)}.\ \square\end{eqnarray*}

 \vspace{4mm}\noindent\textbf{Example: $p$ is the  geometric distribution.}  For $\alpha\leq \d$ consider $f_p(z)=\frac{1-\alpha }{1-\alpha z}$ Then $f_{p^{(r)}}(z)=\frac{1-\alpha r}{1-\alpha rz}$ when $r \leq\d \alpha  $ and 
$$f_q(z)=\frac{1}{2\alpha }(1-\sqrt{1-4\alpha (1-\alpha )z}),\  f_{q^{(\rho)}}(z)=\frac{1-\sqrt{1-4\alpha (1-\alpha )\rho z}}{1-\sqrt{1-4\alpha (1-\alpha )\rho}}.$$ 
with $r=f_q(\rho)$ or $\rho=\frac{r(1-\alpha r)}{1-\alpha}.$ Note here that, with the notation of \eqref{EXP} we have
$$q=s^{4\alpha(1-\alpha)}_{1/2},\  q^{(\rho)}=s^{4\alpha(1-\alpha)\rho}_{1/2}$$

\vspace{4mm}\noindent\textbf{Example: $q$ is  the Sibuya distribution.} For $\d\leq a=\frac{1}{b}<1$ consider $f_p(u)=\frac{u}{1-(1-u)^b}$ and for $0<r<1$ 
$$f_{p^{(r)}}(u)=(1-(1-r)^b)\times \frac{u}{1-(1-ru)^b}.$$ Denote $\rho=1-(1-r)^b.$ Then \eqref{EXP} defines the progeny of $p^{(r)}.$

\section {Acknowledgements} I am deeply  indebted to 'jandri', an anonymous contributor of the Internet site  'les mathematiques.net' where I raised for the first time 
the question of the positivity of the coefficients $(P_n)$ in \eqref{LESP}. He  suggested  the possibility that the $P_n$ are all positive if and only if $1<b\leq 2$, he observed that $P_5<0 $ if $b>3$ and he even pointed out that for $b=2+10^{-9}$ the first negative $P_n$ occurs for $n=45,$ a strong support for the correctness of the statement of Proposition 1.

\section{References}

\vspace{4mm}\noindent \textsc{Athreya, K.B. and  Ney, P.E.} (1972)
{\it Branching Processes.} Springer, New York.

\vspace{4mm}\noindent \textsc{Grey, D.R.} (1975) Two necessary
conditions for embeddability of a Galton-Watson branching process.
{\it Math. Proc. Cambridge Phil. Soc.} {\bf 78}, 339-343.

\vspace{4mm}\noindent \textsc{Harris, T.H.} (1962) {\it The Theory
of  Branching Processes.} Springer, New York.

\vspace{4mm}\noindent \textsc{Kozubowski, T. and Podg\'orski, K.} (2018)  A generalized Sibuya distribution. 
{\it Ann. Inst. Stat. Math.} {\bf 70}, 855-887.

\vspace{4mm}\noindent \textsc{Sibuya, M.} (1979) Generalized hypergeometric, digamma and trigamma distributions.
{\it Ann. Inst. Stat. Math.} {\bf 31}, 373-390.

\vspace{4mm}\noindent \textsc{Toulouse, P.S.} (1999) {\it Th\`emes
de probabilit\'es et statistique.} Dunod, Paris.

\vspace{4mm}\noindent \textsc{Whittaker, E.T. and Watson, G.N.} (1986) {\it  A Course in Modern Analysis.} Cambridge University Press.

\end{document}